\documentclass{birkjour}
 \newtheorem{thm}{Theorem}[section]
 
 \newtheorem{lem}[thm]{Lemma}
 
 \theoremstyle{definition}
 
 \theoremstyle{remark}

 \numberwithin{equation}{section}

\newcommand{\p}{\mathcal P}
\newcommand{\om}{\Omega}
\newcommand{\M}{\mathcal M}
\newcommand{\g}{\mathbf g}
\newcommand{\K}{\mathcal K}

\newcommand{\R}{\mathbb R}

\newcommand{\h}{\mathcal{H}}

\begin{document}

%
%
%
%
%
%
%
%
%

\title[Nodal sets of Laplace eigenfunctions]
 {Nodal sets of Laplace eigenfunctions:\\ estimates of the Hausdorff measure in \\ dimensions
two and three.}

\author[Logunov]{Alexander Logunov}
\address{School of Mathematical Sciences, Tel Aviv University, Tel Aviv 69978, Israel}
\address{Chebyshev Laboratory, St. Petersburg State University, 14th Line V.O., 29B, Saint Petersburg 199178 Russia}
\email{log239@yandex.ru}
\author[Malinnikova]{Eugenia Malinnikova}
\address{Department of Mathematical Sciences,
Norwegian University of Science and Technology
7491, Trondheim, Norway}
\email{eugenia@math.ntnu.no}

\thanks{A.L.  was  supported was supported in part by ERC Advanced Grant 692616 and ISF Grants 382/15 and 1380/13.
Eu. M. was supported by Project 213638 of the Research Council of Norway.}

\dedicatory{In memory of our teacher Victor Petrovich Havin}

\begin{abstract}
 Let $\Delta_M$ be the Laplace operator on a compact $n$-dimensional Riemannian manifold without boundary.
 We study the zero sets of its eigenfunctions $u:\Delta_M u + \lambda u =0$.  In dimension $n=2$ we refine the 
 Donnelly-Fefferman estimate by showing that
$\h^1(\{u=0 \})\le C\lambda^{3/4-\beta}$ for some $\beta \in (0,1/4)$.
 The proof employs the Donnelly-Fefferman estimate  and a  combinatorial argument, which also gives a lower (non-sharp) bound 
in dimension $n=3$:  $\h^2(\{u=0\})\ge c\lambda^\alpha$ for some $\alpha \in (0,1/2)$. The positive constants $c,C$ depend on the manifold, $\alpha$ and $\beta$ are universal. 
\end{abstract}

\maketitle
\section{Introduction} Let $\Delta_M$ be the Laplace operator on a compact $n$-dimensional Riemannian manifold without boundary.
It was conjectured by Yau, see \cite{Y}, that the nodal sets $E_\lambda=\{u_\lambda=0\}$ of  Laplace eigenfunctions $u_\lambda$, $\Delta_M u_\lambda+\lambda u_\lambda=0$  satisfy the following inequality
\[C_1\sqrt{\lambda}\le \h^{n-1}(E_\lambda)\le C_2\sqrt{\lambda}.\]
This conjecture was proved by Donnelly and Fefferman under the assumption that the Riemannian metric is real-analytic (\cite{DF}). The left-hand side estimate is also proved for smooth non-analytic surfaces by Br\"uning (\cite{B}). 

The previous best known estimate from below for a non-analytic manifold in higher dimensions is 
$$ \h^{n-1}(E_\lambda)\ge C\lambda^{(3-n)/4}, $$ 
which gives a constant for $n=3$.   
The two known approaches are: (1)  follow the ideas of Donnelly and Fefferman and find many balls on the wave-scale $\lambda^{-1/2}$ with bounded doubling index, as it is done in \cite{CM} or (2)  use the Green formula $2\int_{E_\lambda}|\nabla_{_M} u_\lambda|= \lambda \int_M |u_\lambda|$ and the estimate  $\frac{\|u_\lambda\|_\infty}{\|u_\lambda\|_1} \leq C  \lambda^{(n-1)/4}$, see \cite{SZ}. The approach in \cite{CM}  also exploits the Sogge-Zelditch estimates of $L^p$-norms of eigenfunctions.
The following upper estimate in dimension two was established by Donnelly and Fefferman, see \cite{DF1},
$$ \h^1(E_\lambda)\le C\lambda^{3/4}. $$
 In this paper we obtain  tiny improvements  to the estimate from below in dimension three and to the estimate from above in dimension two.
We  show that in dimension $2$
\begin{equation} \label{eq:2d}
\h^1(E_\lambda)\le C\lambda^{3/4-\beta},
\end{equation}
 for some $\beta \in (0,1/4)$. It gives a small refinement to the Donnelly-Fefferman estimate.
The proof of \eqref{eq:2d} relies on the results and methods from \cite{DF,DF1}. Roughly speaking,  the Donnelly-Fefferman argument, which gives the estimate with  $\frac{3}{4}$,
 is combined with a combinatorial argument presented below, which gives the $\beta$ improvement.
The same combinatorial argument shows that in dimension $n=3$
\begin{equation} \label{eq:3d}
 \h^2(E_\lambda)\ge C\lambda^\alpha, 
 \end{equation}
 for some $\alpha >0$. As far is we know it gives the first bound that grows to infinity as $\lambda$ increases, but we note that the latter result is not sharp and can be improved up to the bound $c\sqrt{\lambda}\le \h^{n-1}(E_\lambda)$ conjectured by Yau.

This paper is the first part of the work, which consists of three parts.  Polynomial upper estimates for the Hausdorff measure of the nodal sets in higher dimensions are proved in the second part \cite{L1} by a new technique of propagation of smallness.
The lower bound in Yau's conjecture is proved in the third part \cite{L2}
as well as its harmonic counterpart (Nadirashvili's conjecture). 
We remark that the results in \cite{L1,L2} do not give the estimate \eqref{eq:2d} and all three parts can be read independently.


\section{Toolbox}
\subsection{Inequalities for solutions of elliptic equations}

 Let  $(\M, \g)$ be a smooth Riemannian manifold and $\Delta_{\M}$ be the Laplace operator on $\M$, which is defined by the metric $g$. We always assume that the metric is fixed. In the sequel we consider $\M=M\times\R$, where $M$ is a compact manifold with a given metric, on which we study the eigenfunctions, and $\M$ is endowed with the usual metric of the product. Although $\M$ is not compact itself, we will always work on the compact subset $P=M\times[-1,1]$ of $\M$ where all our estimates are uniform.

  A function $h$ on $\M$ is called harmonic if it satisfies the elliptic equation 
\begin{equation}\label{eq:L}
L (h)=div(\sqrt{g} (g^{ij})\nabla h)=0
\end{equation}
 in local coordinates. More precisely, the Laplace operator on $\M$ is given by $\Delta_{\M}(f)=\frac{1}{\sqrt{g}}div(\sqrt{g} (g^{ij}) \nabla f)$.
Harmonic functions satisfy the maximum and minimum principles and the standard  elliptic gradient estimates, see for example \cite[Chapter 3]{GT}. 
 Further, there exists a constant $C$ such that for any geodesic ball $B(x,r)\subset P$
 \begin{equation} \label{eq:h0}
 |\nabla h(x)| \leq \frac C r  \sup\limits_{ B(x,r)} |h|.
\end{equation}
 The Harnack inequality holds: if $h$ satisfies \eqref{eq:L} and $h>0$ in $B(x,r)$, then for any $y \in B(x,\frac 2 3 r)$
 \begin{equation}
  \frac 1 C h(y)< h(x)< C h(y). 
\end{equation}
The following corollary of the Harnack inequality will be also used later.
 If  $h$ satisfies (\ref{eq:L}) and $h(x) \geq 0$ then 
\begin{equation} \label{eq:h2}
 \sup\limits_{B(x,r)} h \geq c \sup\limits_{B(x,\frac 2 3 r)} |h| 
\end{equation}
for some $c=c(\M )>0$  (it follows from the Harnack inequality applied to the function $\sup_{B(x,r)} h-h$).

 \subsection{Estimates on the wavelength scale.}
Now let $(M, g_0)$ be a compact Riemannian manifold. We consider a Laplace eigenfunction $u$  which satisfies $\Delta_M u = - \lambda u$. Adding a new variable, we  consider  the function $h(\xi,t)= u(\xi) e^{\sqrt \lambda t}$ on the product manifold $\M=M \times \mathbb{R}$. The function $h$ appears to be  harmonic on $\M$.  This observation can be used to claim that on the wave-scale   $\lambda^{-1/2}$ the behavior of the Laplace eigenfunctions  reminds that of harmonic functions. This well-known trick was successfully exploited for example in \cite{L, 
 NPS, M}.

Let $\xi$ be an arbitrary point on $M$. Denote by $B(\xi,r)=B_r(\xi)$ the geodesic ball with center at $\xi$ of radius $r$, when the center of the ball is not important we will omit it in the notation and write $B_r$. 

\begin{lem}\label{le:har}  There exist a small number $\varepsilon = \varepsilon (M) >0$ and constants $C_1=C_1(M)$,  $C_2=C_2(M)$  such that for any eigenfunction $u$, $\Delta_M u = - \lambda u$, and any  $r< \varepsilon \lambda^{-1/2}$ the following inequalities hold\\
\begin{equation}\label{eq:e1}
(a)\ \sup_{B_{r}}|u|\le 2\max_{\partial B_{r}}|u|,\quad (b)\
 \sup_{B_{\frac 1 2 r}}|\nabla u|\le C_1\max_{\partial B_{r}}|u| / r.
\end{equation} 
If, in addition, $u(\xi)\geq 0$ and (therefore) $A=\max_{\partial B_{r}(\xi)}u>0$, then 
\begin{equation}\label{l2} 
\sup_{B_{\frac{2}{3} r }(\xi)}|u|\le C_2A.
\end{equation}
 \end{lem}
The inequalities  \eqref{eq:e1} and  \eqref{l2} follow from the standard elliptic estimates, we provide the proofs for the convenience of the reader.
 
We work in local coordinates on $\M= M \times \R$ and  consider the harmonic function $h(\eta,t)= u(\eta) e^{\sqrt \lambda t}$ on $\M$.  The Laplace operator corresponds (locally) to an elliptic operator $L$ defined on a bounded subdomain $\om$ of $\mathbb{R}^{n+1}$, see above. We choose local coordinates such that the distance on the manifold is equivalent to the Euclidean distance (for example by choosing normal coordinates).  Denote by $G_{\om,L}$ the Green function for $L$ on $\om$. By $|x-y|$ we denote the ordinary Euclidean distance between points $x$ and $y$, locally $|x-y|$ is comparable to the distance between the corresponding points in the Riemannian metric on $\M$. We use the following upper estimate (see \cite{KOW}, \cite{JS}, \cite{HHSM2}) of the Green function:  
\[G_{\om,L}(x,y) \leq \frac{C}{|x-y|^{d-2}},\]
where $d=n+1$ is the dimension of $\M$. The constant $C$ depends on the coordinate chart on $\M$, we  consider a finite set of charts that covers $M\times[-1,1]$.

\begin{proof}[Proof of Lemma \ref{le:har}.]
First we suppose that $\sup_{B_r(\xi)}u>0$ and prove that
\begin{equation}\label{eq:e0} 
\sup_{B_{r}(\xi)} u\le 2 \max_{\partial B_{r}(\xi)} u, 
\end{equation}
if $r< \varepsilon(M) \lambda^{-1/2}$, where  $\varepsilon(M)$ is a sufficiently small positive number, which will be chosen later.
 Put $A=\max_{\partial B_{r}(\xi)}u$ and $K= \sup_{B_{r}(\xi)}u $. Let $\xi_0$ be a point in the closed ball $\overline{B_r}$, where $K$ is attained.  

We  consider the cylinder   $Q= B_r \times (-1/ \sqrt \lambda,1/ \sqrt \lambda)$ on $\M$, on which we have $\sup_Q h \leq e K $.
Without loss of generality we assume that $Q \subset \Omega$.
Define the function $w(\eta,t):= A  + e K \lambda t^2$  and note that $w \geq h$ on $\partial Q$ and $|L w| < C_3 K\lambda  $ in $Q$ for some $C_3=C_3(M)$. Now, consider the difference  $v= h-w$. It is non-positive on  $\partial Q$ and satisfies $v(\xi_0,0)=K-A$ and $|Lv| \leq  C_3 K\lambda $ on $Q$.   

We can decompose $v$ into the sum $v=g_1+g_2$, where $g_1$ is a non-positive harmonic function in $Q$ with $g_1|_{\partial Q} = v|_{\partial Q}$ and $g_2(y)=  \int_Q G_{Q,L}(x,y) Lv(x) dx $.   
   Since $Q\subset  \Omega$, the Green function  satisfies 
\[G_{Q,L}(x,y)\leq  G_{\om,L}(x,y) \leq 
C_4/ |x-y|^{d-2}.\] 
Further, for any $y\in Q$ a simple estimate gives
\[\int\limits_{Q} |x-y|^{2-d} dx\leq C_5 \frac{r}{\sqrt{\lambda}}\leq C_5 \frac{\varepsilon}{\lambda}.\] 
Combining the estimates, we get
\[g_2(y)=  \int\limits_{Q}   G_{Q,L}(x,y)Lv(x) dx  \leq C_3C_4 K\lambda  \int\limits_{Q} |x-y|^{2-d} dx \leq C_6  K \varepsilon.\] 
Hence $g_2(\xi_0,0) \leq C_6  K \varepsilon$. The function $g_1$ is non-positive in $Q$ and therefore $K-A= v(\xi_0,0)= g_1(\xi_0,0) +g_2(\xi_0,0) \leq  C_6 K \varepsilon$. Thus $A\geq K(1-C_6  \varepsilon) $ and the right-hand side is greater than $\frac 1 2 K$ if $\varepsilon$ is chosen sufficiently small, \eqref{eq:e0} follows.
 
 The inequality \eqref{eq:e1} (a) follows from \eqref{eq:e0} if one replaces $u$ by $-u$. 
 Finally, the inequalities \eqref{eq:e1} (b)  and \eqref{l2} are obtained combining  \eqref{eq:e1} (a) with \eqref{eq:h0} and \eqref{eq:h2} respectively, where the last two inequalities are applied to the harmonic function   $h$.
\end{proof}
\subsection{Doubling index}\label{ss:di}
Let  $h$ be a harmonic function on $\M$. Locally $h$ can be considered as a solution to the elliptic equation $Lh=0$. We identify $h$ with a function on the cube $\K^d_\rho=[-\rho,\rho]^d\subset\R^d, d=n+1$. We choose local geodesic coordinates, then the metric is locally equivalent to the Euclidean one and $L$ is a small perturbation of the Euclidean Laplace operator. Let $l$ be a positive odd integer such that $l>2\sqrt{d}$, $l=2l_0+1$. For each cube $q$ in $K^d_\rho$ let $lq$ denote the cube obtained from $q$ by the homothety with the center at the center of $q$ and coefficient $l$. Suppose that $lq\subset K^d_\rho$, then we define the doubling index $N(h,q)$ by
\[\int_{lq}|h(x)|^2dx=2^{N(h,q)}\int_{q}|h(x)|^2dx.\] 
Doubling index was used for estimates of the nodal sets in \cite{DF, DF1} and in many subsequent works.  We will need the following properties of the doubling index.
\begin{lem}\label{l:new}
  	(i) ($L^\infty$-estimate) If   a cube $q$ is inscribed in a ball $B$ (and therefore $lq$ contains $2B$), then
\[
\sup_{\frac 4 3 B}|h|\le C_7 2^{ N(h,q)/2}\sup_{B}|h|,\]
for some positive $C_7=C_7(M)$.\\
 (ii) (Monotonicity property) There exists a positive integer $A=A(d)$, a constant $C_0=C_0(d)>1$ and a positive number $\rho=\rho(M)$  such that if $q_1$ and $q$ are cubes that are contained in  $\K^d_\rho$,  and $Aq_1\subset q$ then $N(h, q_1)\le C_0N(h, q)$.\\
 \end{lem}
\begin{proof}
 
(i) 
Indeed, we have 
\[\int_{2 B}|h|^2\le \int_{l q}|h|^2=2^{N(h,q)}\int_{q}|h|^2\le 2^{N(h,q)}\int_{B}|h|^2.\] 
Clearly, $\int_B|h|^2\le(\sup_{B}|h|)^2|B|$. Further, by an elliptic estimate for $h$ , 
$\sup_{\frac 4 3 B}|h| \le C \left(\int_{2B} |h|^2\right)^{1/2}|B|^{-1/2}.$ The inequality follows.\\
(ii) The monotonicity property  is left without a proof. We refer to \cite{GL} and \cite{M}  for the proof of the monotonicity property of the doubling index defined through integrals over concentric geodesic spheres instead of cubes. Using this, it is not difficult to derive the monotonicity property of doubling index for cubes instead of spheres. 
\\
	\end{proof}
   
\section{Inscribed balls and a local estimate of the volume of the nodal set}


The aim of this section is to estimate  from below the volume of the nodal set of an eigenfunction  $u$ of the Laplace operator  in a geodesic ball of radius comparable to the wavelength $ \lambda^{-1/2}$, where  $\Delta_M u+\lambda u=0$.  The estimates presented in this section are very far from being sharp.

 Let us fix a point $O$ on $M$ and assume $u(O)=0$. Denote by $|x|$ the distance from the point $x$ to $O$. We will consider the geodesic ball $B_r$ of radius   $r \leq \varepsilon\lambda^{-1/2}$ and with center at $O$, where $\varepsilon= \varepsilon(M)$ is  chosen so that the inequalities (\ref{eq:e1}) and (\ref{l2}) hold.

\begin{lem}
Assume that $\sup_{B_{\frac r 2 }}|u| \le 2^N\sup_{B_{\frac r 4 }}|u|$, where $N$ is a positive integer, $N\ge 4$. 
Then
\begin{equation} \label{eq:est}
\h^{n-1}\{|x|\le r/2, u(x)=0\}\ge cr^{n-1}N^{2-n},
\end{equation}
for some positive $c=c(M)$.
\end{lem}
\begin{proof} Applying \eqref{l2}, one can deduce 
\[\frac{\max_{\partial B_{r/2}}u}{\max_{\partial B_{3r/8}}u}\le C_2\frac{\sup_{B_{r/2}}|u|}{\sup_{B_{r/4}}|u|} \leq C_2 2^N.\]

Let $S_j=\{x:|x|=r_j=r(\frac 3 8  + \frac j {8N})\}$, $m^+_j=\max_{S_j} u$ and $m^-_j=\min_{S_j} u$, $j=0,1,...,N$.  Recall that $u$ is zero at $O$.  It follows from the weak maximum principle \eqref{eq:e0}  that 
\[m^-_j<0,\ m^+_j>0\quad {\text{ and}}\quad  m^+_{j}\le 2m^+_{j+1},\ |m^-_{j}|\le 2|m^-_{j+1}|.\]
We consider the ratios $\tau_j=m^+_{j+1}/m^+_j$, $j=0,...,N-1$. Then each $\tau_j\ge 1/2$ and 
\[\tau_0...\tau_{N-1}=\frac{\max_{\partial B_{r/2}}u}{\max_{\partial B_{3r/8}}u}\le C_2\frac{\sup_{B_{r/2}}|u|}{\sup_{B_{r/4}}|u|}\le C_22^{N}.\] 
Therefore at most $N/4$ of the ratios $\tau_0,...,\tau_{N-1}$ are greater than $C_3=C_3(C_2)$. Similarly, at most $N/4$ of the ratios  $|m^-_{j+1}|/|m^-_j|$ are greater than $C_3$. Hence there are at least $N/2$ numbers  $k,\ 0\le k\le N-1$ such that $m^+_{k+1}\le C_3m^+_k$ and $|m^-_{k+1}|\le C_3|m^-_k|$. We want to show that for each such $k$ there is a ball of radius $cr/N$ and centered on the sphere $S_k$ where $u$ is positive.
   
Indeed, let $x_0$ be such that $|x_0|=r_k$ and $u(x_0)=m_k^+=\max_{\{|x|=r_k\}}u(x)$ and let $b$ be the ball centered at $x_0$ with radius $\frac{r}{16N}$. Then 
\[\sup_b u \leq  \max_{\{|x|\leq r_{k+1}\}}u(x) \leq C_1 \max_{\{|x|=r_{k+1}\}}u(x) \leq  C_4 m_k^+.\] 
Applying \eqref{l2} we see that $\max_{\frac 1 2 b} |u| \leq C_5 m_k^+$. Taking into account \eqref{eq:e1} (b) and $u(x_0)=m_k^+$, we deduce  that $u$ is positive in a smaller ball of radius $c_1 r /N$ centered at $x_0$. 

 Similarly, we can find a ball of radius $c_1 r /N$ with center on $S_k$ where $u$ is negative.
 Thus the spherical layer $\{x: r_{k-1}<|x|<r_{k+1}\}$ contains two balls of radius $c_1 r /N$ where $u$ has  opposite signs.
Then 
\[\h^{n-1}\{x: r_{k-1}<|x|<r_{k+1}: u(x)=0\}\ge c_2 \left(\frac{r}{N}\right)^{n-1}.\] 
The last inequality holds for at least $N/2$ numbers $k$, so \eqref{eq:est} follows.
\end{proof}


\section{Combinatorial argument}\label{ss:ca}

We need the following lemma about the doubling index that was defined in Section \ref{ss:di}. This lemma holds for an arbitrary function $h\in L^2$, not necessarily harmonic.

\begin{lem} \label{l3}
 Let a cube $Q$ be partitioned into $(Kl)^d$ equal cubes $q_i$ with side length $\frac{1}{Kl}$ (where $l$ is the odd integer from the definition of the doubling index and $K$ is an arbitrary positive integer). Put $N_{\min}= \min\limits_{i} N(h, q_i)$, the minimum is taken over those cubes $q_i$ of the partition for which $lq_i\subset Q$,  and assume that $N_{\min} \geq 2d\ln l/\ln 2$. Then $N(h, \frac{1}{l}Q) \geq  \frac1{2}K N_{\min}$.
\end{lem}

\begin{proof}
Define $Q_j = \frac{K+j(l-1)}{Kl} Q$ for $j=0,1,\dots ,K$, in particular $Q_0=\frac{1}{l}Q$ and $Q_{K}=Q$. 

We know that $\int_{lq_i}|h|^2 \geq 2^{N_{\min}}\int_{q_i}|h|^2$ for each $q_i$ and therefore 
\[2^{N_{\min}}\int_{Q_j}|h|^2\le \sum\limits_{q_i \subset Q_j}\int_{lq_i}|h|^2\le l^d \int_{Q_{j+1}}|h|^2,\]
 since the union of the (open) cubes $lq_i, q_i\subset Q_j$, is contained in  $Q_{j+1}$ and covers each point of $Q_{j+1}$ with multiplicity at most $l^d$. 

Further, the inequality $N_{\min}\ge 2d\ln l/\ln 2$ implies $2^{N_{\min}/2}\int_{Q_j}|h|^2\leq \int_{Q_{j+1}}|h|^2$. Finally, multiplying the last inequalities for $j=0,..., K-1$,  we obtain 
\[\int_{Q_{K}}|h|^2 \geq 2^{KN_{\min}/2} \int_{Q_{0}}|h|^2= 2^{K N_{\min}/2} \int_{\frac 1 l Q}|h|^2.\]

\end{proof}

Suppose now that $h$ is a harmonic function on $\M=M\times\R$.
Given a cube $c$, define $\tilde{N}(h,c)=\sup_{c'\subset c} N(h,c')$, where the supremum is taken over all subcubes $c'$ of the cube $c$.
The monotonicity property implies 
\[ \tilde N\left(h, \frac{1}{A}c\right)\le C_0N(h, c),\]
when $c$ is contained in $\K^{d+1}_\rho$ and $\rho$ is small enough. If a cube $c$ contains a cube $c'$, then  $\tilde N(h, c) \geq \tilde N(h, c')$.

Our aim is to divide the cube $q=\K^{d+1}_{\rho/l}$ into  small cubes and estimate the number of cubes with large doubling constants. 

\begin{lem} \label{c1}
  Let $h$ be a solution to $Lh=0$ in $q$. There exist constants $B_0=B_0(d,L)$ and $\delta=\delta(d)>0$ such that if  the cube $q$ is partitioned into $B>B_0$ equal subcubes, then at least half of these subcubes $c$ satisfy
\[\tilde{N}(h,c)\le\max \left\{ \frac{\tilde{N}(h,q)}{B^{\delta}}, \frac{2d\ln l}{\ln 2} \right\}.\] 
\end{lem}

\begin{proof}
Let $N_0=\tilde{N}(h,q)$. We will do the partition step by step.  On the zero step we have one cube $q$ with $\tilde N(h,q)=N_0$.  We fix $A$ and $C_0$ from the monotonicity property of the doubling index and choose an integer $K$ such that $K>4C_0$. 

 On the first step we divide $q$ into $Y=[lKA]^d$ subcubes. First, divide $q$ into $[lK]^d$ subcubes. By Lemma \ref{l3} at least one subcube $c$ satisfies $N(h,c) \leq 2N_0/K$ if $N_0$ is large enough. Then   $\tilde N(h, \frac{1}{A}c) \leq 2C_0N_0/K \leq N_0/2$. Thus if we divide $q$ into $[lKA]^d$ subcubes, then 
at least one subcube will have $\tilde N \leq N_0/2$ and all other subcubes will have $\tilde N \leq N_0$.

 On the second step we will repeat the partition procedure in each subcube $c$ from the first step. Then at least one subcube $c'$ of $c$ will have $\tilde N(h,c')\leq \tilde N(h,c)/2 $. Also $\tilde N(h,c'')\leq \tilde N(h,c) $ for any other subcube $c''$ of $c$.  

Going from the step with number $j-1$ to the step with number $j$, we take any cube $c$ from the previous step and divide it into $Y$ equal subcubes. In each cube with $\tilde{N}(h,c)\le N_0/2^s$, $1\le s\le j$, we get at least one cube $c'$ with $\tilde{N}(h,c')\le N_0/2^{s+1}$ and for other cubes in $c$ we have $\tilde{N}(h,c')\le N_0/2^s$. 

Using the standard induction  argument, one can see that on the $j$-th step there is one cube with the doubling index less than or equal to $N_0/2^j$, ${j \choose 1}(Y-1)L^d$ other cubes with the indices less than or equal to $N_0/2^{j-1}$, and so on, with  ${j \choose k}(Y-1)^{j-k}L^d$ other cubes with the indices bounded by $N_0/2^{k}$, $k\ge 0$ (assuming that $N_0/2^{j}\geq 2d\ln l/\ln 2$). The sum  $\sum_{k=0}^j {j \choose k}(Y-1)^{j-k}=(1+(Y-1))^j$ is the number of all cubes on the $j$-th step.

 Let $\xi_1,\dots, \xi_j$ be i.i.d. random variables such that $\p(\xi_1 = 1)= 1/Y$ and $\p(\xi_1 = 0)= (Y-1)/Y$.
By the law of large numbers 
\[\p\left( \frac {\sum_{i=1}^j \xi_i }{j}> \frac{1}{2Y}\right) \to 1\ {\text{ as}}\  j \to \infty.\] 
If $j$ is sufficiently large, then  
\[ \frac 1 2\leq  \p\left(\sum_{i=1}^j \xi_i  \geq \frac{j}{2Y}\right) = \sum_{j \geq k \geq \frac{j}{2Y}} \p\left(\sum_{i=1}^j \xi_i = k\right) =  \sum_{ k \geq \frac{j}{2Y}} {j \choose k}\frac{(Y-1)^{j-k}}{Y^j}  .\]

 We conclude that at least  half of all cubes on the $j$-th step have doubling indices  bounded by $N_0/ 2^{\frac{j}{2Y}}  $. Let $B=[lKA]^{jd}=Y^j$, then $N_0/ 2^{\frac{j}{2Y}}  \leq N_0/ B^\delta$, where $\delta=\delta(Y)$ is a positive number such that $Y^{\delta} < 2^{\frac{1}{4Y}}$ and thus $\delta$ depends only on the dimension $d$. Here we have assumed  that $j>j_0$ to apply the law of large numbers and we have also  assumed that $N_0/ B^{\delta} \geq 2d\ln l/\ln 2 $ to apply Lemma \ref{l3}.
 \end{proof}

\section{Estimates of the nodal sets of eigenfunctions}
\subsection{Lower estimate in dimension three}
Suppose now that $u$ is the Laplace eigenfunction, $\Delta_M u+\lambda u=0$ on $M$, where $M$ is a smooth Riemannian three dimensional manifold. Using the standard trick, we consider the manifold $\M=M\times \R$ and a new function $h(\xi,t)=u(\xi)e^{\sqrt{\lambda}t}$, which satisfies $\Delta_\M h=0$. We therefore work on a four-dimensional manifold.

We fix a cube $Q$ on $M$  and consider the cube $\tilde Q = Q \times I$ on $\M$, where $I$ is the interval centered at the origin with the length equal to the side length  of $Q$, we choose $Q$ small enough such that a chart for $Q$ in normal coordinates is contained in some $K^4_\rho$. 

The Donnelly-Fefferman estimate, see \cite{DF}, implies  that $\tilde N(u,Q)\le C\sqrt{\lambda}$ for some $C=C(M)$ if the diameter of $Q$ is less than $c(M)$, and therefore $\tilde N(h,\tilde Q) \leq C_1\sqrt{\lambda}$.
See also \cite{M} for the explanation of the Donnelly-Fefferman estimate via the three sphere theorem for harmonic functions.

 We partition $\tilde Q$ into $B$ smaller cubes $\tilde q$
 with the side length of order $\lambda^{-1/2}$,  such that for each small cube $\tilde q$ there is a zero of $h$ within $\frac{1}{10} \tilde q$ (it is well-known, see \cite{DF}, that the nodal set of $u$ is $c\lambda^{-1/2}$ dense on $M$).   Then $B \sim [c\sqrt{\lambda}]^4 \sim c_1\lambda^2$ and $B$ is large enough when $\lambda>\lambda_0$. 
 
By Lemma \ref{c1}, half of all small cubes have doubling indices bounded by $C \sqrt \lambda / B^\delta  \leq C_1 \lambda^{1/2 - 2\delta}$.
 In each small cube of the wavelength size $C/\sqrt{\lambda}$ the doubling index for $h$  is comparable to the doubling index for the function $u$ on the projection of the cube to $M$, since $h(x,t)=u(x)e^{\sqrt{\lambda}t}$.  Then  at least one half of the small cubes   of size $C/\sqrt \lambda$ in $Q$ have  doubling indices bounded by $C_2 \lambda^{1/2 - 2\delta}$. In each such cube $q$ we can find a smaller subcube $ q' $ with diameter $\frac{\varepsilon}{\sqrt \lambda}$ such that  $u$ is equal to $0$ at the center of $ q'$. Then combining Lemma \ref{l:new} (i) and  the estimate \eqref{eq:est}, we obtain 
\[\h^2(\{u=0 \}\cap q') \geq  \frac{c_2}{\lambda N(u, q')} \geq c_3 \lambda^{-3/2} \lambda^{2\delta}.\]
The number of such cubes is comparable to $ \lambda^{3/2}$. Thus $\h^2(\{u=0\}) \geq c_4 \lambda^{2\delta}$.
 
\subsection{Upper estimate in dimension two} Following \cite{DF1}, using local isothermal coordinates in a geodesic disk of radius $r$, we transform the eigenfunction $u$, $\Delta_M u+\lambda u=0$ to a function $f$ in the unit ball of $\R^2$ that satisfies $\Delta_0 f+\lambda r^{2}\psi f=0 $, where $\Delta_0$ is the Euclidian Laplacian and $\psi$ is a bounded function (the bound depends on the metric).  

We will combine the combinatorial argument from Section \ref{ss:ca} with the following estimate for the length of the nodal set by Donnelly and Fefferman, \cite{DF1}. Let $Q$ be the unit square.

{\it Suppose that $g:Q\to \R$ satisfies $\tilde{N}(g, Q)\le \Gamma $, $\Delta g=\Gamma\psi g$, where $\psi$ is a  function in $Q$ with sufficiently small $L^\infty$-norm. Then \[\h^1(x\in \frac{1}{100}Q: g(x)=0)\le C\Gamma.\] }

In \cite{DF1} this estimate is  applied on the scale $\lambda^{-1/4}$: for any square $q$ on $M$ with side $\sim \lambda^{-1/4}$ one can consider a function  $u(\lambda^{-1/4} x)$ and apply  the estimate with $\Gamma \sim \lambda^{1/2}$ (using  that the doubling index for any cube is bounded by $C\lambda^{1/2}$) to see that $H^{n-1}(\{u=0\}\cap q) \leq C \lambda^{1/4}$. Summing the estimates over such  cubes covering $M$, one has $H^1(\{u=0\})\le C\lambda^{3/4}$.

 However a combinatorial argument will show that  very few cubes with side $\sim \lambda^{-1/4}$ have doubling indices comparable to $ \lambda^{1/2}$, in fact, most of the cubes have significantly smaller doubling indices. We are going to refine the global length estimate via combining the combinatorial argument 
and the Donnelly-Fefferman estimate on various scales. 
\begin{lem} Fix  a geodesic ball $B$ on the surface with isothermal coordinates and let $q$ be a square in $B$ with side-length $\sim \lambda^{-1/4}$. Then
\begin{equation}\label{str}
\h^1(\{u=0\} \cap q)\le C \tilde{N}(u, 100 q)^{1/2}.
\end{equation}
\end{lem}
 \begin{proof} Denote $\tilde{N}(u,100 q)$ by $N_0$. Let us 
divide $q$  into squares with side-length $\sim N_0^{1/2}\lambda^{-1/2}$. In each of those the doubling index is bounded by $N_0$ and rescaling such small squares to unit squares and applying the estimate of Donnelly and Fefferman with $\Gamma= N_0$, we bound the length of the nodal set in such small square by $CN_0^{3/2}\lambda^{-1/2}$. The number of such squares is $\sim \frac{\lambda^{1/2}}{N_0} $. Then 
$\h^1(\{u=0\} \cap q)\le CN_0^{3/2}\lambda^{-1/2}\lambda^{1/2}N_0^{-1}=CN_0^{1/2}$
\end{proof}

Now let $\K=\K^2_\rho$ be a square such that $100\K$ lies in (the chart for) $B$, the side-length of $\K$ depends only on the geometry of the surface $M$ and does not depend on $\lambda$. 
We partition $\K$ into squares with side-length $\lambda^{-1/4}$, then for each such square $q$ we have $\h^1(\{u=0\} \cap q)\le C \tilde{N}^{1/2}(u, 100 q_0) $ and summing up over all squares $q$ in the partition of $\K$, we obtain
$$\h^1(\{u=0\}\cap \K)\le C \sum_{q\subset \K}\tilde{N}^{1/2}(u, 100 q).$$

Further,  we consider the harmonic  extension  $h(t,x)=e^{\sqrt{\lambda}t}u(x)$ of $u$ and let $\tilde\K=\K^3_\rho=\K\times[-\rho,\rho]$.  Note that $N(h,\tilde q)\ge N(u,q)$, whenever $q$ is the projection of $\tilde q$ onto $M$ and then the same inequality holds for $\tilde{N}$.  

 Let $Y$ be a sufficiently large integer defined in Section \ref{ss:ca}. Choose an integer $j$ such that $Y^j \sim \lambda^{3/4}$.  We partition $\tilde \K$  into $Y^j \sim \lambda^{3/4}$ subcubes with side-length $\sim \lambda^{-1/4}$.  According to Section \ref{ss:ca} these cubes can be divided into $j$ groups $G_0,...,G_j$ such that  $\tilde{N}(h,\tilde q)\le N_02^{s-j}$ for each cube $\tilde q\in G_s$, where $N_0:=\tilde{N}(h, K_0)\le C\sqrt{\lambda}$ and the number of cubes in $G_s$ is ${j \choose s} (Y-1)^s$.
However  we need to replace  $\tilde{N}(h,\tilde q)$ by $\tilde{N}(h,100\tilde q)$ in the estimate for a number a cubes  in order to estimate the sum  $\sum_{q\subset \K}\tilde{N}^{1/2}(u, 100 q)$. It can be done by changing the parameter $l$ in the definition of the doubling index in Section \ref{ss:di}. The doubling index with a parameter $l$ in a cube $100c$ can be estimated by  the doubling index with a parameter $10000l$ in a cube $c$.
We therefore have $\tilde{N}(h,100\tilde q)\le C  N_02^{s-j}$ for each cube $q\in G_s$, here we abuse the notation $\tilde{N}$ for a doubling index with the modified $l$ and denote it by the same letter.

Finally, we apply the inequality $\tilde N(u,q)\le \tilde N(h, \tilde q)$, where $q$ is the projection of $\tilde q$, and estimate $\tilde N^{1/2}(u, 100 q)$ by the average of the corresponding quantities over $Y^{j/3}$ cubes $\tilde q$ with the projection $q$. We obtain
\[
\h^1(\{u=0\}\cap \K)\le C\sum_{q\subset \K}\tilde{N}^{1/2}(u, 100 q)\le CY^{-j/3}\sum_{\tilde q\subset \tilde \K}\tilde{N}^{1/2}(h, 100\tilde q).\]
Further we partition all cubes $\tilde q$ into the groups $G_s$,
\begin{multline*}\sum_{\tilde q\subset \tilde \K}\tilde{N}^{1/2}(h, 100\tilde q)=\sum_{s=0}^j\sum_{\tilde q\in G_s}\tilde{N}^{1/2}(h, 100\tilde q)\\
\le C\lambda^{1/4}\sum_{s=0}^j{j\choose s}(Y-1)^s2^{(-1/2)(j-s)}=C\lambda^{1/4}(Y-1+2^{-1/2})^j.\end{multline*}
We have $Y^j=c\lambda^{3/4}$, then $Y-1+2^{-1/2}=Y^{1-\eta}$ for some $\eta=\eta(Y)>0$ and $\h^1(\{u=0\})\le C_M\lambda^{3/4(1-\eta)}$.

\subsection*{Acknowledgment}

We are grateful to Lev Buhovsky and Mikhail Sodin, who   
read the first
draft of this paper and made very helpful suggestions and comments.

This work was started when the first author visited NTNU and the second author visited  the Chebyshev Laboratory (SPBSU). The work was finished at TAU and Purdue University. We are grateful to these institutions for their hospitality and for great working conditions.

\end{document}